%begin macros
\parskip=8pt
\font\eightrm=cmr8  
\font\eighttt=cmtt8
\magnification=\magstephalf
\parindent=0pt
\overfullrule=0in
%end macros
\bf
\noindent
HOW JOE GILLIS DISCOVERED COMBINATORIAL SPECIAL FUNCTION THEORY
\smallskip
\noindent
\it
Doron Zeilberger
\smallskip
\it
\qquad \qquad \qquad \qquad \qquad
\qquad Girsa d'yankuta [la mishtakcha]\footnote{$^1$}
{\eightrm Aramaic: The lesson of infancy is not forgotten.} \qquad 
(Talmud, Tractate Shabbat 22b)
\medskip
\rm
{\bf Nov. 19, 1993:} When I checked my mail this morning, I was
shocked to learn that Joe Gillis died last night in his sleep.
He was eighty two years old.
 
Only last May, during my last visit to Israel, he discussed with me his
research plans. He was hoping to generalize results on Hausdorff dimension
and what are now called fractal sets, that he had found back in 
the mid-thirties.
 
Joe had a great influence on my mathematical development, as he
did on many generations of Israeli mathematicians, who got their first
taste of advanced mathematics through {\it Gilyonot Lematematika},
an outstanding mathematical magazine, in Hebrew, `` for high-school
students and amateurs''. This magazine, 
which Joe edited for many years, had
an extensive problem section. The problems periodically challenged
the best and brightest among us .
Joe also initiated the Israeli Math
Olympiad, and was the coach of the Israeli I.M.O. team for a very long time.
 
However, Joe's influence to mathematics at large
was far greater than that, in particular, with respect to
my own two specialties: Combinatorics and Special Functions. Gillis,
in his seminal paper with Even[2],
initiated the extremely fruitful
marriage of these hitherto unrelated subjects, from which
was soon born
the flourishing new field of combinatorial special functions (e.g.[3][5]).
I feel that the story of how this unison came 
to be, narrated to
me, years ago, by Joe himself, must be recorded for posterity, since it
testifies not only to Joe's genius, but to the genius of the human
spirit.
 
{\bf 1928:} Joe spent his last year in high-school preparing for the
competition for the
coveted scholarship to Trinity College,
Cambridge, that he subsequently won. The textbook he studied 
was Chrystal's famous {\it Algebra}[1]. One of the problems
discussed there particularly appealed to Joe. It was the
classical {\it derangements} problem:
{\it In how many ways can one stuff $n$ different letters, 
in the corresponding $n$ envelopes, in such
a way that no letter gets sent to the right address?}. The well-known
answer, given in Chrystal's text, is that this number, $D(n)$,
equals $n!(1/0!-1/1!+1/2!-1/3! + \dots + (-1)^n/n!) =[n!/e]$.
 
Joe, being the bright youngster that he was, started to wonder
what happens if there are {\it multiple} letters addressed to each
address. In other words, what can one say about $D( n_1 , \dots , n_k)$,
the number of ways of stuffing $n_1$ letters addressed to $1$, \dots,
$n_k$ letters addressed to $k$, into the corresponding 
$n_1 + \dots + n_k$ envelopes, in such a way that no letter gets
to the right destination? Of course, he realized that
 
$$
D(n_1,n_2)=  \delta_{n_1 , n_2} \quad ,
\eqno(1)
$$
 
but was unable to find a ``closed form'' expression, even for the case $k=3$.
Failing this, he went on to establish {\it recurrence relations},
that enabled him to compile a table for $D(n_1 , n_2 , n_3)$, 
for small (and not so small) values  of the arguments, starting from
the obvious ``initial conditions'' (1).
Having accomplished this, and realizing that there
probably is no reasonable formula for $D(n_1, n_2 , n_3)$, he went
on to ``bigger and better''\footnote{$^2$}
{\eightrm After completing his
undergraduate studies with distinction, he went on to write
a brilliant dissertation under Besicovitch; was one of the first
collaborators of Erd\"os; was stationed in Bletchley Park; made important
contributions to fluid dynamics; and so on..., but this is a different story.}
things, or so it seemed then.
 
{\bf 1960:} About one-third-century later, 
and long after he ``changed fields'' to
``applied'' mathematics (spurred initially by his desire to contribute
to the welfare of the then young state of Israel, to where he immigrated in
the late forties), he encountered a ``practical'' problem. In the course of
trying to solve a certain differential equation, he needed to compute
the following ``linearization coefficients'' for the Laguerre polynomials:

$$
E( n_1 , n_2 , n_3 ) =
(-1)^{(n_1 + n_2 + n_3)} \int L_{n_1}(x) L_{n_2}(x) L_{n_3}(x) e^{-x} dx \quad.
\eqno(2)
$$
 
Once again, he was unable to find a ``closed form'' expression.
However, he and George Weiss[4] obtained recurrence relations for
the $E(n_1 , n_2 , n_3)$ which enable one to compile
a table of these for any specified range of the arguments, 
obviating the need to integrate every time anew.
 
{\bf 1975:} A decade and a half later, as he was browsing through his
old paper[4], Joe made a {\it connection.}
He had seen these recurrence relations for $(2)$, 
established $15$ years earlier,
{\it much  much} earlier! They were identical (up to some trivial
change of notation) to the recurrences he established for 
$D(n_1 , n_2 , n_3 )$,
almost half-a-century before, during his last year of high school!
Matching the obvious initial conditions at $n_1=0, n_2=0, n_3=0$, for
which $E$ coincides with $(1)$(due to the orthogonality of the Laguerre
polynomials), it followed that ([2])
 
$$
D( n_1 , n_2 , n_3 ) =  E( n_1 , n_2 , n_3 )  \quad.
\eqno(3)
$$
 
Thus began the beautiful field of combinatorial special function theory
(e.g. [3][5].)\footnote{$^3$}
{\eightrm
The referee pointed out that the connection between combinatorics and
function theory goes back to Euler, Gauss, and Jacobi. However the connection
between combinatorics and the classical special functions of mathematical 
physics, was first made, as far as I am aware of, by Joe Gillis.
I wish to thank the referee for many valuable comments.}
\medskip
\noindent
{\bf References}
 
{\bf 1. } G. Chrystal, {``Algebra''}, vol. {\bf 2}, reprinted by
Chelsea, N.Y., N.Y.
 
{\bf 2.} S. Even and J. Gillis, {\it Derangements and Laguerre polynomials},
Proc. Cambridge Phil. Soc. {\bf 79}(1976), 135-143.
 
{\bf 3.} D. Foata, {\it Combinatoire de identiti\'es sur le polyn\^omes
orthogonaux}, Proc. Inter. Congress of Math.[Warsaw. Aug. 16-24, 1983],
Warsaw, 1983.
 
{\bf 4.} J. Gillis and G. Weiss, {\it Products of Laguerre polynomials},
M.T.A.C. (now Math. Comp.), {\bf 14}(1960), 60-63.
 
{\bf 5.} J. Zeng,
{\it  Weighted derangements and the linearization coefficients of orthogonal
Sheffer polynomials}, Proc. London Math. Soc. (3) {\bf 65}(1992), 1-22.

{\obeylines
Dept. of Mathematics, 
Temple Univ., 
Philadelphia, PA 19122, USA. 
{\eighttt zeilberg@math.temple.edu .}
} 
\bye